\documentclass[12pt,a4paper]{article}
\usepackage[utf8]{inputenc}
\usepackage[T1]{fontenc}
\usepackage{amsmath}
\usepackage{amsfonts}
\usepackage{amssymb,amsthm}
\usepackage{graphicx}
\usepackage{mathrsfs}
\usepackage{tikz}
\usepackage{hyperref}
\usepackage{color}
\usepackage{etoolbox}
\usepackage{appendix}
\usepackage{enumitem}
\usepackage{comment}
\newtheorem{thm}{Theorem}[section]

\newtheorem{pro}[thm]{Proposition}
\newtheorem{cor}[thm]{Corollary}
\newtheorem{lem}[thm]{Lemma}
\newtheorem{rem}[thm]{Remark}
\newcommand{\R}{\mathbb{R}}

\newcommand{\na}{\nabla}

\newcommand{\Om}{\Omega}
\newcommand{\Omb}{\overline{\Om}}
\newcommand{\p}{\partial}
\newcommand{\s}{\sigma}

\makeatletter
\@addtoreset{equation}{section} 
\patchcmd{\maketitle}{\@fnsymbol}{\@alph}{}{}  
\makeatother
\title{An interpolation approach to $L^{\infty}$ a priori estimates for elliptic problems with nonlinearity on the boundary }

\makeatletter
\newcounter{author}
\renewcommand*\author[1]{%
\stepcounter{author}%
\ifnum\c@author=1
\gdef\@author{#1}%
\else
\xdef\@author{\unexpanded\expandafter{\@author\and#1}}%
\fi
\csgdef{author@\the\c@author}{#1}}
\newcommand*\email[1]{%
\csgdef{email@\the\c@author}{#1}}
\newcommand*\address[1]{%
\csgdef{address@\the\c@author}{#1}}
\AtEndDocument{%
\xdef\author@count{\the\c@author}%
\c@author=1
\print@authors}
\newcommand*\print@authors{%
\ifnum\c@author>\author@count
\else
\print@author{\the\c@author}%
\advance\c@author by 1
\expandafter\print@authors
\fi}
\newcommand*\print@author[1]{%
\par\medskip
\begin{tabular}{@{}l@{}}%
\textsc{Addresses of \csuse{author@#1}}\\
\csuse{address@#1}\\
\textit{E-mail address}:
\href{b:mailto:\csuse{email@#1}}{\csuse{email@#1}}
\end{tabular}}
\makeatother
 {\small
 \author{Maya Chhetri}
 \address{
 Department of Mathematics and Statistics\\
 University of North Carolina  Greensboro \\
 Greensboro, NC, USA}
 \email{m\_chhetr@uncg.edu}
 \author{Nsoki Mavinga}
 \address{
 Department of Mathematics and Statistics\\
 Swarthmore College, Swarthmore, PA, USA}
 \email{nmaving1@swarthmore.edu}
 \author{Rosa Pardo}
 \address{Department of Mathematical Analysis and Applied Mathematics\\
 Universidad Complutense de Madrid, Madrid, Spain}
 \email{rpardo@ucm.es}
 }

\date{\today}

\begin{document}
\maketitle
{ }\makeatletter{\renewcommand*{\@makefnmark}{}
} 

{ }\makeatletter{\renewcommand*{\@makefnmark}{}
}

\begin{abstract}
We establish  an explicit $L^\infty(\Om)$ a priori estimate for weak solutions to subcritical 
elliptic problems with nonlinearity on the boundary, in terms of the powers of their $H^1(\Om)$ norms.
To prove our result, we combine in a novel way Moser type estimates 
together with elliptic regularity and Gagliardo--Nirenberg interpolation inequality. We illustrate our result with an application to subcritical problems satisfying Ambrosetti-Rabinowitz condition.
\end{abstract}

\section{Introduction}
In this work, our goal is to establish an explicit $L^\infty(\Om)$ a priori estimate  of weak solutions to elliptic equations with nonlinear boundary conditions of the form  
\begin{equation}
\label{eq:pde}
\left\{
\begin{array}{rcll}
-\Delta u +u &=&  0 \quad &\mbox{in}\quad \Om\,;\\
\frac{\p u}{\p \eta} &=&  f(x, u)\quad &\mbox{on}\quad \p \Om\,.
\end{array}
\right. 
\end{equation}
Here $\Om \subset \R^{N}\ (N >2)$  is a bounded domain with Lipschitz boundary $\p \Om$, and  $\p/\p\eta
:=\eta(x)\cdot\nabla$ denotes the outer normal derivative on the boundary $\p\Om$.  The nonlinearity $f :\p\Om\times \R \to \R$ is a Carathéodory function, that is,  $f(\cdot,s)$ is measurable for all $s\in\R$,  and $f(x,\cdot)\ $  
 is   continuous for a.e. $x \in \p\Om$. Further, for a.e. $x \in \p\Om$ and for all $s \in \R$, the nonlinearity  $f$ satisfies the subcritical growth condition: 
\begin{equation}
\label{growth:f:subcritical}
|f\,( x,s)|\le B_0(1+|s|^{\, p})\,, \quad 1 < p < 2_*-1
\end{equation}
for some positive constant $B_0$, where $2_*:=\frac{2(N-1)}{N-2}$ is the critical exponent for the trace operator. 
\par 
\par We say that a weak solution $u$ of \eqref{eq:pde} has $L^\infty(\Om)$ {\it a priori estimate} if  $u \in L^\infty(\Om)$ and $\|u\|_{L^\infty(\Om)} \leq M$, where the constant $M=M(u,f,\Om)$. 
\par  It was established in  \cite[Thm.~3.1]{Marino-Winkert_2019} that any weak solution of a general quasilinear problem which includes \eqref{eq:pde}, belongs to $L^{\infty}(\Om) \cap L^{\infty}(\p\Om)$. It turns out that if $u$ is a weak solution of \eqref{eq:pde} with the property that $u \in L^{\infty}(\Om)$, then $u \in C(\overline{\Omega})$ as well (see Proposition \ref{pro:reg}).

Our goal here is to describe the bound $M(u,f,\Om)$ explicitly in terms of powers of the $H^1(\Om)$ norm of $u$.
To this end, we establish the following result.

\begin{thm}
\label{thm:estim} 
Suppose $f$ satisfies \eqref{growth:f:subcritical}. Then there exists $C_0>0$, depending only on  $N$ and $\Om$, such that for every weak  solution  $u$ of \eqref{eq:pde}
\begin{equation}\label{asymp:estim:scalar}
\|u\|_{L^{\infty}(\Om)} 
\le C_0\Big(1+\|u\|_{H^1(\Om)}\Big)^{A}\,,
\end{equation}
where 
\begin{equation}\label{def:A}
A:=\frac{2_*-2}{(2_*-1)-p}=\frac{2}{N-p(N-2)}>0.
\end{equation}
\end{thm} 
For the  semilinear and quasilinear cases with the nonlinearities  in the interior satisfying the zero Dirichlet boundary condition, see \cite{Pardo_JFPTA}  and \cite{Pardo_RACSAM} respectively, for similar explicit $L^{\infty}(\Om)$  estimates. To the best of our knowledge, the explicit estimate \eqref{asymp:estim:scalar} for elliptic problems with nonlinearity on the boundary has not been studied before. In general, problems with nonlinear boundary conditions are not widely studied. As a result, one  of the difficulties of dealing with such problems  depend on the obtainment of the appropriate regularity results. 

Uniform $L^\infty(\Om)$ a priori bounds (that is, when $M$ is independent of the solution $u$) are generally studied for positive solutions, see \cite[Theorem 3.7]{Bon-Ros_2001}. They used the blow-up method introduced by \cite{Gidas-Spruck} with the aid of the nonexistence results of \cite{Hu} to nonlinear boundary conditions.
The availability of a uniform $H^1(\Om)$ a priori bound result for sign-changing weak solutions of \eqref{eq:pde} will complement the result of Theorem~\ref{thm:estim} resulting in a uniform $L^\infty(\Om)$ a priori bound. This is an open problem.

As an illustration we establish that, for nonlinearities satisfying an Ambrosetti - Rabinowitz condition (see \cite{Ambrosetti_Rabinowitz}), finite energy solutions are uniformly  bounded in  their $H^1(\Om)$-norms. Hence, finite energy solutions will be uniformly $L^\infty(\Om)$ a priori bounded and vice-versa, see Theorem \ref{th:energy}.

\begin{rem}
{\rm 
The proof of Theorem \ref{thm:estim}  also provides $L^{\infty}(\Om)$ estimates of solutions in terms of $L^{2^*}(\Om)$ and $L^{2_*}(\p\Om)$ norms, where $2^*:=\frac{2N}{N-2}$. Specifically, 
it follows from the estimate \eqref{A1:A2} that 
\begin{equation*}
\|u\|_{L^{\infty}(\Om)} \leq C
\bigg(1+\|u\|_{L^{2_*}(\p\Om)}^{\hat{A}_1}\bigg)
\bigg(1+\|u\|_{L^{2^*}(\Om)}^{\hat{A}_2}\bigg) \,,
\end{equation*}
where 
the exponents $\hat{A}_1$ and 
$\hat{A}_2$, given by \eqref{exp}, satisfy $ \hat{A}_1+  \hat{A}_2=A$,
and $C>0$ is independent of $u$,  depending only on  $N$   and $\Om$. 
}
\end{rem}

As a consequence of Theorem~\ref{thm:estim} we have the following two corollaries regarding uniform bound and  convergence to zero with respect to various norms.

\begin{cor}
\label{cor:equiv}
Suppose $f$ satisfies \eqref{growth:f:subcritical}. 
Let $\{u_n\}$ be a sequence of weak solutions  to \eqref{eq:pde}.
Then the following four statements are equivalent:
 \begin{itemize}
 \item[]\rm(i)~  \ $\|u_n\|_{L^{2_*}(\p\Om)}\le C_1$;
\qquad \rm (ii)~\ $\|u_n\|_{H^1(\Om)}\le C_2$;
 \item[]\rm (iii)~\ $\|u_n\|_{L^{\infty}(\Om)}\le C_3$;
\qquad\ \rm (iv)~\ $\|u_n\|_{C^{\mu}(\overline{\Om})}\le C_4$. 
 \end{itemize}
for some constants $C_i$ independent of $n$, $i=1,\cdots 4$.
\end{cor}

\begin{cor}
    \label{cor:eq}
Suppose $f$ satisfies \eqref{growth:f:subcritical}.     
Let $\{u_n\}$ be a sequence of weak solutions  to \eqref{eq:pde}. Then the following four statements are equivalent as $n\to\infty$:
\begin{itemize}
 \item[]\rm(i)~  \ $\|u_n\|_{L^{2_*}(\p\Om)}\to 0$;
\qquad \rm (ii)~\ $\|u_n\|_{H^1(\Om)}\to 0$;
 \item[]\rm (iii)~\ $\|u_n\|_{L^{\infty}(\Om)}\to 0$;
\qquad\ \rm (iv)~\ $\|u_n\|_{C(\overline{\Om})}\to 0$. 
 \end{itemize}

\end{cor}
We 
combine  Moser type estimates (\cite{Brezis_Kato, Moser, DKN_book, Marino-Winkert_2019, Struwe})
together with elliptic regularity and Gagliardo--Nirenberg interpolation inequality (\cite{Brezis_2011, Nirenberg}) in a novel way 
to prove Theorem~\ref{thm:estim}.

The paper is organized in the following way. In Section \ref{sec:prelim} we collect some known results.
Section \ref{sec:proof:thm:estim:cor} contains the proof of Theorem~\ref{thm:estim}. Section \ref{sec:cor} is devoted to the proofs of Corollary~\ref{cor:equiv} Corollary~\ref{cor:eq}. Finally, in Section \ref{sec:app} we give an application   of our results to finite energy solutions, see Theorem~\ref{th:energy}. 

Throughout this paper $C$ will denote a constant which may depend on $f$, $\Om$ and $N$, but it is independent of the solution $u$, and may vary from line to line.

\section{Preliminaries}
\label{sec:prelim}
In this section, we define weak solutions, recall the regularity of weak solutions to the linear problem with non-homogeneous Neumann boundary conditions, and state a result that guarantees in particular $L^{\infty}$ regularity of weak solutions.
\medskip

A function $u \in H^1(\Om)$ is a {\it weak solution} to \eqref{eq:pde} if
\begin{equation}
\label{def:weak:scalar}
\int_{\Omega} \left[\nabla u \nabla \psi + u \psi \right]= \int_{\p\Om} f(x, u)\psi , \qquad \forall \psi \in H^1(\Om)\,.
\end{equation}
The left hand side of \eqref{def:weak:scalar} is well defined since $u, \psi \in H^1(\Om)$. Also, by the continuity of the trace operator $u, \psi \in L^{2_*}(\p\Om)$.  
The Carathéodory assumption and the growth condition  \eqref{growth:f:subcritical} imply that $f(\cdot, u(\cdot)) \in L^{\frac{2_*}{p}}(\p\Om)$. Using H\"{o}lder's inequality, we get  
\begin{equation}\label{f:weak}
    \int_{\p\Om} |f(x, u(x))\psi|\le \|f(\cdot, u(\cdot))\|_{L^{(2_*)'}(\p\Om)}\|\psi\|_{L^{2_*}(\p\Om)},
    \end{equation}
where $(2_*)':=\frac{2_*}{2_*-1}$ is the conjugate of $2_*$. Then the right hand side of \eqref{f:weak} is finite since  $p < 2_*-1$ and   $(2_*)'=(\frac{2_*}{p})(\frac{p}{2_*-1})<\frac{2_*}{p}$,  and the integrals in \eqref{def:weak:scalar} are well defined.
 
Next, we
consider the linear problem 
\begin{equation}
\label{lbp} 
\left\{ \begin{array}{rcll}
-\Delta  v +v &=&0 & \qquad \mbox{in } \Om\,;  \\
\frac{\p   v}{\p \eta}&=&  \mathit{h} &  \qquad \mbox{on } \p
\Om\,,
\end{array}\right.
\end{equation}
where $h \in L^q(\p\Omega)$ for $q \geq 1$. We denote the solution operator  by  
$T: L^{q}(\p\Om)\to W^{1,m}(\Om) \mbox{ with }  \mathit{T}\,\mathit{h}:=v\,$ and 
\begin{equation}
\label{ineq:m:h}
\|v\|_{W^{1,m}(\Om)}\le  C\|\mathit{h}\|_{L^{q}(\p\Om)},\qquad \text{where }\quad 1 \leq  m\le Nq/(N-1)\,,
\end{equation}
see, for instance \cite{Ama_1976, Lad-Ura_1968,  Mavinga-Pardo_PRSE_2017, Necas1986} for more details. 

It is known that the trace operator
$
\Gamma : W^{1,m}(\Om ) \rightarrow L^{r}(\p
\Om ) 
$
is a continuous linear operator for every $r$ satisfying 
$
\frac{N-1}{r}\ge \frac{N}{m}-1$, and compact  if $\frac{N-1}{r}> \frac{N}{m}-1$, see \cite[Ch.~6]{KufnerJohnFucik77}.
\smallskip
Now, we define  the resolvent operator 
$S: L^q(\p\Om) \to L^r(\p\Om)$  by $S\mathit{h}:=  \Gamma \big(\mathit{T}\,\mathit{h}\big)=\Gamma v\,,$
given schematically by
$$
L^q(\p\Om)\stackrel{T}{\longrightarrow} W^{1,m}(\Om)\stackrel{\Gamma}{\longrightarrow}L^r(\p\Om)
$$
for any $q \geq 1$ and for all $r$ satisfying $\frac{N-1}{r}\ge \frac{N-m}{m}$ with $1 \leq m\le Nq/(N-1)$. 
Note that if $\frac{N-1}{r}> \frac{N-m}{m}$ then $S$ is compact by the compactness of $\Gamma$.
\medskip

The following result states the regularity of the solution to the linear problem \eqref{lbp}. 

\begin{lem}
\label{lem:reg:lbp}
Let $N\ge  2$ and $\mathit{h}\in L^q (\p \Om )$
with $q\ge  1$. Then, the unique solution $v=\mathit{T}\mathit{h}$ of the linear problem \eqref{lbp} satisfies the following:
\begin{enumerate}
\item [(i)] If $1\le  q<N-1$, then $\Gamma v\in L^{r}(\p\Om)$
for all $1\le  r\le  \frac{(N-1)q}{N-1-q}$ and the map $
S:L^q(\p\Om)\to L^r(\p\Om)$ is
continuous for $1\le  r\le  \frac{(N-1)q}{N-1-q}$  and compact for 
$1\le  r<\frac{(N-1)q}{N-1-q}$. 

\item [(ii)] If $q=N-1$, then $\Gamma v\in L^{r}(\p\Om)$ for 
all $r\ge  1$ 
and the map $ S:L^q(\p\Om)\to L^r(\p\Om)$ is
continuous and compact for $1\le  r<\infty$. 

\item [(iii)] If $q>N-1$, then $v\in C^{\mu}(\overline{\Om})$ with
$\|v\|_{C^{\mu}(\overline{\Om})}\le 
C\|h\|_{L^q(\p\Om)}$ for 
some $\mu\in (0,1)$. Moreover, 
$\Gamma v\in C^{\mu}(\p\Om)$ 
and the map $S:L^q(\p\Om)\to C^{\mu}(\p\Om)$ is  continuous and compact.

\item [(iv)] If $\mathit{h}\in C^{\mu}(\p\Om)$, then $v\in C^{2,\alpha }(\Om)\cap C^{1,\alpha }(\overline{\Om})$. 
\end{enumerate}
\end{lem}
\begin{proof}
See \cite[Lemma 2.1]{Arrieta-Pardo-RBernal_2007} for proofs of (i)-(iii). See  \cite{Lad-Ura_1968, Necas, Necas1986} for the proof of (iv).
\end{proof}

Part of the following regularity result was established in \cite[Thm.~3.1]{Marino-Winkert_2019} using Moser’s iteration technique. This enables us to estimate $L^{\infty}(\Om)$ norm of  weak solutions, which is generally only available to classical solutions. 
\begin{pro}
\label{pro:reg} 
Suppose $f:\p\Om \times \R\to\R$ is a Carath\'eodory function such that there exists a constant $B_1$ satisfying
\begin{equation}
\label{f:growth:2}
|f(x, s)|\le B_1(1+|s|^{2_*-1})
\end{equation}
for a.e. $x \in \p\Om$ and for all $s \in \R$. If $u$ is a weak solution of  \eqref{eq:pde},   then  $u\in C^{\mu}(\overline{\Om})\cap W^{1,m}(\Om)$ for any $0<\nu<1$ and  $1<m<\infty$. Moreover,
\begin{equation}\label{infty:cont}
 \|u\|_{L^{\infty}(\p\Om)}\le \|u\|_{C(\Omb)}=\|u\|_{L^{\infty}(\Om)}.    
\end{equation}
\end{pro} 
\begin{proof}
 It follows from  \cite[Thm.~3.1]{Marino-Winkert_2019} that $u \in L^{\infty}(\Om) \cap L^{\infty}(\p\Om)$. Then using the elliptic regularity result we get that  $u\in C^{\mu}(\overline{\Om})\cap W^{1,m}(\Om)$ for any $0<\mu<1$ and  $1<m<\infty$. The inequality \eqref{infty:cont} then follows from \cite[P.~83]{KufnerJohnFucik_book}).
\end{proof}

\section{Proof of Theorem~\ref{thm:estim}}
\label{sec:proof:thm:estim:cor}
Let $u $ be a   weak solution of \eqref{eq:pde}. By  Proposition~\ref{pro:reg},  $u \in C(\Omb)$. 
Fix $q$ such that
\begin{equation}\label{q}
  q>\max\left\{N-1, \frac{2_*}{p}\right\}. 
\end{equation} 
Then, using \eqref{growth:f:subcritical} and \eqref{infty:cont}, we get
\begin{align}
\label{asymp:ineq:q:infty}
\|f(u)\|_{L^q(\p\Om)}^q &\le C
\int_{\p\Om}  \big(1+|u|^{p}\big)^q\le 
C \int_{\p\Om}\big(1+|u|^{pq}\big)  \\
&= C \int_{\p\Om}\Big(1+|u|^{2_*}|u|^{pq-2_*}\Big) 
\le  C\bigg(1+\|u\|_{L^{\infty}(\Om)}^{pq-2_*}\|u\|_{L^{2_*}(\p\Om)}^{2_*}\bigg)\,.
\nonumber
\end{align}
Therefore,  \eqref{ineq:m:h} and \eqref{asymp:ineq:q:infty} imply
\begin{equation}
\label{asymp:ineq:infty:2_*}
\|u\|_{W^{1,m}(\Om)}\le C\, \bigg(1+\|u\|_{L^{\infty}(\Om)}^{pq-2_*}\|u\|_{L^{2_*}(\p\Om)}^{2_*}\bigg)^{1/q}\,,
\end{equation}
taking $m:=\frac{Nq}{N-1}>N$ see \eqref{q}.  Using 
the Gagliardo-Nirenberg interpolation inequality, see \cite{Nirenberg} or \cite[3.C~Ex.~3, p.~313-314]{Brezis_2011}, 
there exists $C(\Om, N, m)$ such that 
\begin{equation}
\label{asymp:ineq:gagliardo}
\|u\|_{L^{\infty}(\Om)} \leq C\|u\|_{W^{1,m}(\Om)}^{\sigma} \|u\|_{L^{2^*}(\Om)}^{1-\sigma}
\end{equation}
where

\begin{equation}\label{sigma}
   \frac{1}{\sigma}:=1+2^{*} \left( \frac{1}{N} - \frac{1}{m}\right)\,. 
\end{equation}
Now substituting \eqref{asymp:ineq:infty:2_*} into \eqref{asymp:ineq:gagliardo}, we get	
	\begin{align}
	\label{asymp:ineq:infty:gagliardo}
	\|u\|_{L^{\infty}(\Om)} &\leq C
	\bigg(1+\|u\|_{L^{\infty}(\Om)}^{pq-2_*}
	\|u\|_{L^{2_*}(\p\Om)}^{2_*}\bigg)^{\s/q}
	\|u\|_{L^{2^*}(\Om)}^{1-\sigma} \nonumber\\
	&\leq C
	\bigg(1+\|u\|_{L^{\infty}(\Om)}^{\sigma(p-\frac{2_*}{q})}
	\|u\|_{L^{2_*}(\p\Om)}^{2_*\s/q}\bigg)
	\|u\|_{L^{2^*}(\Om)}^{1-\sigma}\,.
\end{align}

To establish the estimate, we consider two cases:  $\|u\|_{L^{\infty}(\Om)} > 1$ and   $\|u\|_{L^{\infty}(\Om)} \le 1$.
\par If $\|u\|_{L^{\infty}(\Om)} > 1$,  dividing both sides of \eqref{asymp:ineq:infty:gagliardo} by $\|u\|_{L^{\infty}(\Om)} ^{\sigma(p-\frac{2_*}{q})}$, 
we get
\begin{align}
	\|u\|_{L^{\infty}(\Om)}^{1-\sigma(p-\frac{2_*}{q})} &\leq C
	\bigg(\frac1{\|u\|_{L^{\infty}(\Om)} ^{\sigma(p-\frac{2_*}{q})}}+
	\|u\|_{L^{2_*}(\p\Om)}^{2_*\s/q}\bigg)
	\|u\|_{L^{2^*}(\Om)}^{1-\sigma}\nonumber \\
	\label{for:remark:scalar}
	&\leq C
	\bigg(1+
	\|u\|_{L^{2_*}(\p\Om)}^{2_*\s/q}\bigg)	\,
\|u\|_{L^{2^*}(\Om)}^{1-\sigma} \,,
\end{align}
since $p>2_*/q$ by \eqref{q}. For $m=\frac{Nq}{N-1}$, we have
\begin{equation}
\label{m:q}
\frac{2^*}{m}
=\frac{2N}{N-2}\frac{N-1}{Nq}=\frac{2(N-1)}{N-2}\frac{1}{q}=\frac{2_*}{q}\,.
\end{equation}
Then using  \eqref{sigma} and \eqref{m:q}, we deduce
\begin{align*}
	\frac{1}{\sigma}+\frac{2_*}{q}&=1+2^{*}\left( \frac{1}{N} - \frac{1}{m}\right)+\frac{2_*}{q}
	=1+\frac{2}{N-2}=\frac{N}{N-2}	
	\,,
\end{align*}
hence,
\begin{equation}\label{denom}
1>    1-\sigma\left(p-\frac{2_*}{q}\right)=\s\left(\frac{N}{N-2}-p\right)>0
\end{equation}
since 
$ p < 2_{*}-1=\frac{N}{N-2}.$ Then it follows from \eqref{for:remark:scalar} that
\begin{align}
\label{A1:s:A2}
\|u\|_{L^{\infty}(\Om)} &\leq C
\bigg(1+\|u\|_{L^{2_*}(\p\Om)}^{\hat{A}_1}\bigg)
\|u\|_{L^{2^*}(\Om)}^{\hat{A}_2}\\
&\leq C
\bigg(1+\|u\|_{L^{2_*}(\p\Om)}^{\hat{A}_1}\bigg)
\bigg(1+\|u\|_{L^{2^*}(\Om)}^{\hat{A}_2}\bigg) \,,
\nonumber
\end{align}
where $\hat{A}_1$ and $\hat{A}_2$ are positive exponents given by
\begin{equation}\label{exp}
\hat{A}_1:=\frac{2_*\s/q}{\s\left(\frac{N}{N-2}-p\right)}=\frac{\frac{2^*}{m}}{\frac{N}{N-2}-p},\
\hat{A}_2:=\frac{1-\s}{\s\left(\frac{N}{N-2}-p\right)}=\frac{\frac{2^*}{N}-\frac{2^*}{m}}{\frac{N}{N-2}-p}.
\end{equation}

On the other hand, if $\|u\|_{L^{\infty}(\Om)} \le 1$, then \eqref{asymp:ineq:infty:gagliardo} yields, 
\begin{align}
	\label{asymp:ineq:infty:gagliardo:2}
	\|u\|_{L^{\infty}(\Om)} &\leq C
	\bigg(1+\|u\|_{L^{2_*}(\p\Om)}^{2_*\s/q}\bigg)
	\,\|u\|_{L^{2^*}(\Om)}^{1-\sigma}\\
 \label{asymp:ineq:infty:gagliardo:3a}
 &\leq C
	\bigg(1+\|u\|_{L^{2_*}(\p\Om)}^{{\hat{A}_1}}\bigg)
	\,\|u\|_{L^{2^*}(\Om)}^{1-\sigma}\\
 &\leq C
	\bigg(1+\|u\|_{L^{2_*}(\p\Om)}^{{\hat{A}_1}}\bigg)
	\,\bigg(1+\|u\|_{L^{2^*}(\Om)}^{1-\sigma}\bigg)\nonumber\\
 \label{A1:A2}
 &\leq C
\bigg(1+\|u\|_{L^{2_*}(\p\Om)}^{\hat{A}_1}\bigg)
\bigg(1+\|u\|_{L^{2^*}(\Om)}^{\hat{A}_2}\bigg) \,,
\end{align}
where the exponents $\hat{A}_1$ and $\hat{A}_2$ are as given in \eqref{exp}. In the inequalities \eqref{asymp:ineq:infty:gagliardo:3a} and \eqref{A1:A2}, we used the fact that for any $a\le b$, there exists a constant $c>0$ such that
$1+x^a\le c(1+ x^{b})$,  
for all $x\ge 0$. We employed \eqref{denom}, and let  $x=\|u\|_{L^{2_*}(\p\Om)}$, $a=2_*\s/q$, $b=\hat{A}_1>2_*\s/q$ to get \eqref{asymp:ineq:infty:gagliardo:3a},
and let   
$x=\|u\|_{L^{2^*}(\Om)}$, $a=1-\s$ and $b=\hat{A}_2>1-\s$ to get \eqref{A1:A2}.

Then,  Sobolev embedding and the continuity of the trace operator yields
\begin{align*}
	\|u\|_{L^{\infty}(\Om)} &
	\le C\, \Big(1+\|u\|_{H^1(\Om)}\Big)^{A}\,,
\end{align*}
with 
$$
	A:=\hat{A}_1+\hat{A}_2=\frac{\frac{2^*}{N}}{\frac{N}{N-2}-p}=\frac{2}{N-(N-2)p} >0\,,     
$$
as desired.
This completes the proof. 
\hfill $\Box$
\section{Proofs of Corollary~\ref{cor:equiv} and Corollary~\ref{cor:eq}}
\label{sec:cor}
\begin{proof}[Proof of Corollary \ref{cor:equiv}]
We  prove that (i)$\implies$(ii)$\implies$(iii)$\implies$(iv)$\implies$(i).

Assume that there exists a constant $C_1>0$  such that $\|u_n\|_{L^{2_*}(\p\Om)}\le C_1$, where $C_1$  is independent of $n$.
Then by definition of weak solution, see \eqref{def:weak:scalar}, Hölder's inequality and \eqref{growth:f:subcritical},
\begin{align}
\|u_n\|_{H^{1}(\Om)}^2 &=\int_{\Om}\left[ |\na u_n|^2 + u_n^2 \right]=\int_{\p\Om} f(x, u_n)u_n \nonumber\\
&\le \|f(\cdot, u_n(\cdot))\|_{L^{(2_*)'}(\p\Om)}\,\|u_n\|_{L^{2_*}(\p\Om)}   
\label{est:h1:2_*}\\
&\le \tilde{C}\left(|\p\Om| + \int_{\p\Om}|u_n|^{2_*} \right)^{1-\frac{1}{2_*}}\,\|u_n\|_{L^{2_*}(\p\Om)} \le (C_2)^2,\nonumber
\end{align}
where $C_2$ is independent of $n$, and  (ii) holds.

Then, (iii) follows since the estimate  \eqref{asymp:estim:scalar} of Theorem~\ref{thm:estim} yields 
$\|u_n\|_{L^{\infty}(\Om)} \leq C_3
$ for some  $C_3>0$  independent of $n$. Moreover, by \eqref{infty:cont}, 
$f(\cdot,u_n(\cdot))\in L^q(\partial \Om)$, for any $q> N-1$. 
By Lemma \ref{lem:reg:lbp} part (iii), $u_n\in C^{\mu}(\overline{\Om})$  for any $\mu \in (0, 1)$, and 
$\|u_n\|_{C^{\mu}(\overline{\Om})}\le   C\|f(\cdot,u_n(\cdot))\|_{L^q(\partial \Om)}$.
Using \eqref{growth:f:subcritical}   we get  $\|u_n\|_{C^{\mu}(\overline{\Om})}\le C_4$, completing (iv). 
If (iv) holds then clearly (i) holds with some constant independent of $n$.
\end{proof}

\begin{proof}[Proof of Corollary \ref{cor:eq}]
We  prove that (i)$\implies$(ii)$\implies$(iii)$\implies$(iv)$\implies$(i).
If (i) holds, that is, $\|u_n\|_{L^{2_*}(\p\Om)}\to 0,$ then (ii) follows from \eqref{est:h1:2_*}.
If (ii) holds, it follows from continuous Sobolev embedding that $\|u_n\|_{L^{2^*}(\Om)}\to 0$. Using \eqref{denom} we get $1- \s<\hat{A}_2$, where  $\hat{A}_2$ and $ 1- \s$ 
are given by \eqref{exp} and \eqref{sigma} respectively. Then  the estimates 
\eqref{A1:s:A2} and  \eqref{asymp:ineq:infty:gagliardo:3a} 
imply that there exists $C>0$, independent of $u$,  
such that 
\begin{align}
\label{conv:infty}
	\|u_n\|_{L^{\infty}(\Om)} &\leq C
	\bigg(1+\|u_n\|_{L^{2_*}(\p\Om)}^{\hat{A}_1}\bigg)
	\,\|u_n\|_{L^{2^*}(\Om)}^{1- \s} \to 0,
\end{align}
where the positive exponent $\hat{A}_1$ is given by \eqref{exp},  resulting in (iii). If (iii) holds, by \eqref{infty:cont} $\|u_n\|_{C(\Omb)}=\|u_n\|_{L^{\infty}(\Om)} \to 0$ and (iv) follows. If (iv) holds, \eqref{infty:cont} again implies 
 $ \|u_n\|_{L^{\infty}(\p\Om)}\le \|u_n\|_{C(\Omb)} \to 0
 $, hence  $ \|u_n\|_{L^{2_*}(\p\Om)}\to 0$ and thus (i) holds. This completes the proof.
\end{proof}

\section{An application}
\label{sec:app}
In this section, we discuss an application of our results to finite energy solutions.

We say that a sequence $\{u_n\}\subset H^1(\Omega )$ of   weak solutions  to \eqref{eq:pde} has {\it uniformly bounded energy} if there exists a constant $c_o>0,$ such that $J[u_n]\le c_o,$
where $J$ is the associated energy functional defined by
\begin{equation*}
J[u]  := \frac{1}{2}\int_{\Om}\Big( |\nabla u|^2 +u^2\Big)-
\int_{\p\Om} F(x,u)\,=\frac{1}{2}\|u_n\|^2_{H^1(\Omega)} -
\int_{\p\Om} F(x,u)\,,
\end{equation*}
with $F(x,t):=\int_{0}^{t}f(x,s)\, ds$.

We say that the {\it Ambrosetti--Rabinowitz condition} holds if there exists  two constants $\theta>2,$ and $s_0>0$ such that
\begin{equation}\label{AR}
\tag{AR}
\theta F(x,s)\le sf(x,s),\qquad \forall x\in\Om,\quad \forall |s|>s_0.
\end{equation}

Assuming that Ambrosetti--Rabinowitz condition holds, the next Theorem states that a sequence of solutions to \eqref{eq:pde} is uniformly $L^\infty(\Om)$  {\it a priori} bounded if and only if  it has uniformly bounded energy.

\begin{thm}
\label{th:energy}
Let $\{u_n\}\subset H^1(\Omega )$ be a sequence  of  weak solutions  to \eqref{eq:pde}.
Assume that $f:\p\Om \times\R\to \R$  is a  Carathéodory function  satisfying  \eqref{growth:f:subcritical}  and  \eqref{AR}.
\smallskip
\par Then there exists a uniform constant $C>0$ (depending  only on   $\Omega,$ $N$ and $f$, independent of $u$) such that 
\begin{equation}
\label{uniform:bound}
\|u_n\|_{L^\infty (\Om)}\leq C,
\end{equation}
if and only if $\{u_n\}\subset H^1(\Omega )$  has uniformly bounded energy.	
\end{thm}

\begin{proof}[Proof of Theorem \ref{th:energy}] Let $\{u_n\}\subset H^1(\Omega )$ be a sequence  of  weak solutions  to \eqref{eq:pde} satisfying \eqref{uniform:bound}. 
By \eqref{growth:f:subcritical}, $|F(x,u_n)|\le B(1+|u_n|^{p+1})$ with $p+1<2_*$. Therefore, 
\begin{align*}\label{var:pb:C2}
\big|J[u_n]\big|&\le\frac{1}{2}\int_{\Om}\Big( |\nabla u|^2 +u^2\Big) +
\int_{\p\Om} \big|F(x,u_n)\big|\, dx\\
&\le \frac{1}{2}\|u_n\|^2_{H^1(\Omega)}+ B\|u_n\|^2_{L^{2_*}(\partial\Omega)} +B'\,,
\end{align*}
where $B'>0$ depends only on $f$ and $|\p\Om|$.
Corollary \ref{cor:equiv} implies that $\{u_n\}$ is uniformly bounded in both  $ H^1(\Omega)$-norm and $ L^{2_*}(\p\Omega )$-norm and hence $J[u_n] \leq c_o$ for some $c_o>0$ independent of $n$.
\par Now, suppose that there exists a positive constant $c_o>0$ such that  
\begin{equation*}\label{var:pb:C1}
J[u_n]
=\frac{1}{2}\|u_n\|^2_{H^1(\Omega)} -
\int_{\p\Om} F(x,u_n)\, dx\le c_o.
\end{equation*}
Observe that the condition  \eqref{AR} implies that there exists a constant $C>0$ such that 
\begin{equation*}\label{F:uf:2}
\int_{\p\Om} F(x,u_n)\, dx\le \frac{1}{\theta }\int_{\p\Om} u_nf(x,u_n)\, dx +C.
\end{equation*}
Hence,
\begin{equation*}\label{uf}
\frac{1}{2}\|u_n\|^2_{H^1(\Omega)}-
\frac{1}{\theta } \int_{\p\Om} u_nf(x,u_n)\, dx\le C.
\end{equation*}
By taking $u_n$ as a test function in \eqref{def:weak:scalar} and combining with the previous estimate, we get 
\begin{equation*}\label{uf:2}
\left(\frac{1}{2}-\frac{1}{\theta }\right)\|u_n\|^2_{H^1(\Omega)}\le C\,.
\end{equation*}
Then, using Theorem~\ref{thm:estim}, we conclude that $\|u_n\|_{L^{\infty}(\Om)} \leq C'$, where $C'$ independent of $n$. This completes the  the proof.
\end{proof}

\section*{Acknowledgements}
This material is based upon work supported by the National Science Foundation under Grant No. 1440140, while the authors were in residence at the Mathematical Sciences Research Institute in Berkeley, California, during the month of June of 2022.
\par The first author  was  supported by a grant from the Simons Foundation 965180. The second author was supported by James A. Michener Faculty Fellowship.
The third author is supported by grants PID2019-103860GB-I00, and PID2022-137074NB-I00,  MICINN,  Spain, and by UCM, Spain,  Grupo 920894.

\end{document}